\documentclass[11pt, a4paper]{article}
\usepackage{amsmath}
\usepackage{wasysym}
\usepackage{latexsym}
\usepackage{amsfonts}
\usepackage{mathrsfs}
\usepackage{amssymb}
\usepackage{ifsym}
\usepackage{dsfont}
\usepackage[all]{xy}
\usepackage{authblk}

\newtheorem{Definitions1}{Definition}[section]

\newtheorem{Theorems1}{Theorem}[section]

\newtheorem{Coroll1}[Theorems1]{Corollary}
\newtheorem{Lemma1}[Theorems1]{Lemma}

\newtheorem{Quest1}{Question}[section]

\newenvironment{proof}[1][Proof]{\begin{trivlist}
\item[\hskip \labelsep {\bfseries #1}]}{\end{trivlist}}

\newcommand{\qed}{\nobreak \ifvmode \relax \else
      \ifdim\lastskip<1.5em \hskip-\lastskip
      \hskip1.5em plus0em minus0.5em \fi \nobreak
      \vrule height0.75em width0.5em depth0.25em\fi}

\begin{document}
\title{On the strength of a weak variant of the Axiom of Counting} 
\author{Zachiri McKenzie}
\affil{Department of Philosophy, Linguistics and Theory of Science\\
Gothenburg University\\
\texttt{zachiri.mckenzie@gu.se}}
\maketitle

\begin{abstract}
In this paper $\mathrm{NFU}^{-\mathrm{AC}}$ is used to denote Ronald Jensen's modification of Quine's `New Foundations' Set Theory ($\mathrm{NF}$) fortified with a type-level pairing function but without the Axiom of Choice. The axiom $\mathrm{AxCount}_\geq$ is the variant of the Axiom of Counting which asserts that no finite set is smaller than its own set of singletons. This paper shows that $\mathrm{NFU}^{-\mathrm{AC}}+\mathrm{AxCount}_\geq$ proves the consistency of the Simple Theory of Types with Infinity ($\mathrm{TSTI}$). This result implies that $\mathrm{NF}+\mathrm{AxCount}_\geq$ proves that consistency of $\mathrm{TSTI}$, and that $\mathrm{NFU}^{-\mathrm{AC}}+\mathrm{AxCount}_\geq$ proves the consistency of $\mathrm{NFU}^{-\mathrm{AC}}$.  
\end{abstract}

\section[Introduction]{Introduction}

The Axiom of Counting ($\mathrm{AxCount}$) was introduced by J. Barkley Rosser in \cite{ros78} and asserts that every finite set has the same cardinality as its own set of singletons. When added to Quine's `New Foundations' Set Theory ($\mathrm{NF}$) or Ronald Jensen's variant of $\mathrm{NF}$ that allows urelements ($\mathrm{NFU}$), this axiom proves the comprehension scheme for formulae which fail to be stratified only by virtue of the fact that well-defined types can not be assigned to variables which range only over natural numbers. In the nineteen seventies two natural variants of $\mathrm{AxCount}$ emerged: $\mathrm{AxCount}_\leq$ and $\mathrm{AxCount}_\geq$. $\mathrm{AxCount}_\leq$ asserts that a finite set is no bigger than its own set of singletons, while $\mathrm{AxCount}_\geq$ asserts that a finite set is no smaller than its own set of singletons. It quickly became apparent that many of the strong consequences of $\mathrm{AxCount}$ (over both $\mathrm{NF}$ and $\mathrm{NFU}$) also follow from $\mathrm{AxCount}_\leq$ \cite{hin75, fh09}. In contrast $\mathrm{AxCount}_\geq$ appears to be a much weaker assumption. This paper investigates the strength of $\mathrm{AxCount}_\geq$ over a theory, $\mathrm{NFU}^{-\mathrm{AC}}$, which can be viewed as both a subtheory of $\mathrm{NFU}$ (by which we mean Jensen's system \cite{jen69} supplemented with both with an axiom asserting the existence of a type-level pairing function and the Axiom of Choice) and $\mathrm{NF}$, and shows that this axiom proves the consistency of the Simple Theory of Types with Infinity ($\mathrm{TSTI}$).\\ 
\\
\indent In the context of $\mathrm{NF}$ very little is known about the relative strengths of the theories obtained by adding $\mathrm{AxCount}$, $\mathrm{AxCount}_\leq$ and $\mathrm{AxCount}_\geq$. Steven Orey \cite{ore64} shows that $\mathrm{NF}+\mathrm{AxCount}$ proves the consistency of $\mathrm{NF}$. In \cite{hin75} Roland Hinnion develops techniques that yield lower bounds on the consistency strengths of the theories $\mathrm{NF}+\mathrm{AxCount}$, $\mathrm{NF}+\mathrm{AxCount}_\leq$ and $\mathrm{NF}$ relative to subsystems of $\mathrm{ZFC}$. This paper provides a new lower bound on the consistency strength of $\mathrm{NF}+\mathrm{AxCount}_\geq$ relative to a well-understood $\mathrm{ZF}$-style theory. This lower bound is stronger than any known lower bound on the consistency strength of $\mathrm{NF}$.\\
\\  
\indent There is much clearer picture of the relative strengths of the theories obtained by adding $\mathrm{AxCount}$, $\mathrm{AxCount}_\leq$ and $\mathrm{AxCount}_\geq$ to $\mathrm{NFU}$, largely thanks to Ronald Jensen's consistency proof of $\mathrm{NFU}$ \cite{jen69}. This consistency proof yields the exact strength of $\mathrm{NFU}$ relative to a subsystem of $\mathrm{ZFC}$, and Robert Solovay and Randall Holmes (unpublished) have also computed the exact strength of the theory $\mathrm{NFU}+\mathrm{AxCount}$ relative to a subsystem of $\mathrm{ZFC}$. The paper \cite{mck15} separates the consistency strengths of the theories $\mathrm{NFU}+\mathrm{AxCount}$, $\mathrm{NFU}+\mathrm{AxCount}_\leq$ and $\mathrm{NFU}+\mathrm{AxCount}_\geq$ by showing that $\mathrm{NFU}+\mathrm{AxCount}$ proves the consistency of $\mathrm{NFU}+\mathrm{AxCount}_\leq$, and $\mathrm{NFU}+\mathrm{AxCount}_\leq$ proves the consistency of $\mathrm{NFU}+\mathrm{AxCount}_\geq$. Here it is shown that $\mathrm{NFU}+\mathrm{AxCount}_\geq$ proves the consistency of $\mathrm{NFU}$ answering a question raised in \cite{mck15}.\\       
\\
\indent This paper is only concerned with variants of Quine's `New Foundations' that are fortified with the Axiom of Infinity. It is interesting to note, however, that in \cite{ena06} Ali Enayat investigates the strengths of extensions of the theory $\mathrm{NFU}^{-\infty}$ that is obtained by adding the negation of the Axiom of Infinity to Jensen's modification of $\mathrm{NF}$. Robert Solovay has shown (unpublished) that $\mathrm{NFU}^{-\infty}$ is equiconsistent with the subsystem of arithmetic $I\Delta_0+\mathrm{exp}$. Both $\mathrm{AxCount}$ and $\mathrm{AxCount}_\leq$ are inconsistent with $\mathrm{NFU}^{-\infty}$. The theory $\mathrm{NFU}^{-\infty}+\mathrm{AxCount}_\geq$ is equivalent to the theory $\mathrm{NFUA}^{-\infty}$ which Enayat \cite{ena06} shows is equiconsistent with Peano Arithmetic ($\mathrm{PA}$).

\section[Background]{Background} \label{Sec:Background}

In this section we present the axioms of the theories $\mathrm{NFU}^{-\mathrm{AC}}$ and the Simple Theory of Types ($\mathrm{TSTI}$), and the extensions and subsystems of these theories that we will refer to in the next section of the paper. We also present some necessary facts related to the development of mathematics in $\mathrm{NFU}^{-\mathrm{AC}}$ and outline Roland Hinnion's technique of interpreting well-founded set theories in the set of equivalence classes of topped well-founded extensional relations in $\mathrm{NFU}^{-\mathrm{AC}}$. A detailed development of mathematics in an extension of $\mathrm{NFU}^{-\mathrm{AC}}$ can be found in the textbook \cite{hol98}. We also refer the reader to the monograph \cite{for95} for a treatment of advanced topics in the study of stratified set theories including extensions of $\mathrm{NFU}^{-\mathrm{AC}}$ and $\mathrm{TSTI}$.\\        
\\
\indent Throughout this paper we will use $\mathcal{L}$ to denote the language of set theory: first-order logic endowed with the binary relation $\in$. We use $\mathcal{L}_{\mathrm{PA}}$ to denote the language of arithmetic: first-order logic endowed with binary function symbols $+$ and $\cdot$, and constant symbols $\mathbf{0}$ and $\mathbf{1}$. As usual we will write $\mathrm{PA}$ for the $\mathcal{L}_{\mathrm{PA}}$-theory that consists of all of the axioms of Peano Arithmetic. If $\mathcal{L}^\prime$ is a recursive language and $T$ is a recursively axiomatised $\mathcal{L}^\prime$-theory then we write $\mathrm{Con}(T)$ for the $\mathcal{L}_{\mathrm{PA}}$-formula which asserts that $T$ is consistent.\\
\\
\indent The Simple Theory of Types is the simplification of the underlying system of \cite{rw08} independently discovered by Frank Ramsey and Leon Chistwick. Following \cite{mat01} we use $\mathrm{TSTI}$ to denote the Simple Theory of Types fortified with the Axiom of Infinity. This theory is naturally axiomatised in the many-sorted language with sorts for each $n \in \mathbb{N}$.

\begin{Definitions1}
We use $\mathcal{L}_{\mathrm{TST}}$ to denote the $\mathbb{N}$-sorted language endowed with binary relation symbols $\in_n$ for each sort $n \in \mathbb{N}$. There are variables $x^n, y^n, z^n, \ldots$ for each sort $n \in \mathbb{N}$ and well-formed $\mathcal{L}_{\mathrm{TST}}$-formulae are built-up inductively from atomic formulae in the form $x^n \in_n y^{n+1}$ and $x^n = y^n$ using the connectives and quantifiers of first-order logic.  
\end{Definitions1}

An $\mathcal{L}_{\mathrm{TST}}$-structure $\mathcal{M}$ consists of a function $M$ with domain $\mathbb{N}$ where $M(0), M(1), \ldots$ are the domains of the sorts, and a function $\in^\mathcal{M}$ with domain $\mathbb{N}$ such that for all $n \in \mathbb{N}$, $\in^\mathcal{M}(n) \subseteq M(n) \times M(n+1)$; we write $\mathcal{M}= \langle M, \in^\mathcal{M} \rangle$.

\begin{Definitions1}
We use $\mathrm{TSTI}$ to denote the $\mathcal{L}_{\mathrm{TST}}$-theory with axioms
\begin{itemize}
\item[](Extensionality) for all $n \in \mathbb{N}$, 
$$\forall x^{n+1} \forall y^{n+1} (x^{n+1}= y^{n+1} \iff \forall z^n(z^{n} \in_n x^{n+1} \iff z^n \in_n y^{n+1})),$$  
\item[](Comprehension) for all $n \in \mathbb{N}$ and for all well-formed $\mathcal{L}_{\mathrm{TST}}$-formulae $\phi(x^n, \vec{z})$,
$$\forall \vec{z} \exists y^{n+1} \forall x^n (x^n \in_n y^{n+1} \iff \phi(x^n, \vec{z}))$$
\item[](Infinity) $\exists x^1 \exists f^3(f^3: x^1 \longrightarrow x^1 \textrm{ is injective but not surjective})$.  
\end{itemize} 
\end{Definitions1}

If $n$ is a natural number with $n \geq 4$ then we use $\mathrm{TSTI}_n$ to denote the weakening of $\mathrm{TSTI}$ which only allows formulae to refer to objects with type $<n$.

\begin{Definitions1}
Let $n \in \mathbb{N}$ with $n \geq 4$. We use $\mathcal{L}_n$ to denote the $n$-sorted language endowed with binary relation symbols $\in_k$ for each sort $k < n-1$. There are variables $x^k, y^k, z^k, \ldots$ for each sort $k < n$ and well-formed $\mathcal{L}_n$-formulae are built-up inductively from atomic formulae in the form $x^k \in_k y^{k+1}$ where $k < n-1$, and $x^k = y^k$ where $k < n$, using the connectives and quantifiers of first-order logic.
\end{Definitions1}

If $n \in \mathbb{N}$ then we use $[n]$ to denote the set $\{i \in \mathbb{N} \mid i \leq n\}$ and we use $(n)$ to denote the set $\{i \in \mathbb{N} \mid i < n\}$. An $\mathcal{L}_n$-structure $\mathcal{M}$ consists of a function $M$ with domain $[n-1]$ where $M(0), \ldots, M(n-1)$ are the domains of the sorts, and a function $\in^\mathcal{M}$ with domain $(n-1)$ such that for all $k < n-1$, $\in^\mathcal{M}(k) \subseteq M(k) \times M(k+1)$; we write $\mathcal{M}= \langle M, \in^\mathcal{M} \rangle$.

\begin{Definitions1}
Let $n \in \mathbb{N}$ with $n \geq 4$. We use $\mathrm{TSTI}_n$ to denote the $\mathcal{L}_n$-theory with axioms
\begin{itemize}
\item[](Extensionality) for all $k < n-1$, 
$$\forall x^{k+1} \forall y^{k+1} (x^{k+1}= y^{k+1} \iff \forall z^k(z^{k} \in_k x^{k+1} \iff z^k \in y^{k+1})),$$  
\item[](Comprehension) for all $k < n-1$ and for all well-formed $\mathcal{L}_n$-formulae $\phi(x^k, \vec{z})$,
$$\forall \vec{z} \exists y^{k+1} \forall x^k (x^k \in_k y^{k+1} \iff \phi(x^k, \vec{z}))$$
\item[](Infinity) $\exists x^1 \exists f^3(f^3: x^1 \longrightarrow x^1 \textrm{ is injective but not surjective})$.  
\end{itemize} 
\end{Definitions1}     

If $\sigma_0, \ldots, \sigma_m$ is a proof of a contradiction from $\mathrm{TSTI}$ then there exists an $n \in \mathbb{N}$ such that $\sigma_0, \ldots, \sigma_m$ is a proof of a contradiction from $\mathrm{TSTI}_n$. Combining this with the observation that there exists a binary Turing machine which on input $n \geq 4$ and $k$ decides whether the sentence with G\"{o}del number $k$ is an axiom of $\mathrm{TSTI}_n$ yields 
\begin{equation}
\mathrm{PA} \vdash ((\forall n \geq 4) \mathrm{Con}(\mathrm{TSTI}_n) \iff \mathrm{Con}(\mathrm{TSTI})) 
\end{equation}     

\indent A formula in the language of set theory is said to be stratified if the formula can be turned into a well-formed formula of $\mathcal{L}_{\mathrm{TST}}$ by decorating variables and instances of $\in$ appearing in the formula with types. In 1937 Willard van Orman Quine proposed an axiomatisation of set theory, now dubbed `New Foundations' ($\mathrm{NF}$) after the title of \cite{qui37}, that appears to avoid the set theoretic paradoxes by restricting Cantor's unrestricted comprehension scheme to stratified formulae. In \cite{jen69} Ronald Jensen considers a weakening of $\mathrm{NF}$ that permits both sets and non-sets (urelements) in the domain of discourse. Jensen was able to show that this modification of $\mathrm{NF}$ is consistent relative to a weak subsystem of $\mathrm{ZFC}$ and, unlike $\mathrm{NF}$ (see \cite{spe53}), is consistent with both the Axiom of Choice and the negation of the Axiom of Infinity. In this paper we will stick to the convention of using $\mathrm{NFU}$ to denote Jensen's theory fortified by asserting the existence of a type-level pairing function, and by the Axiom of Choice. We will use $\mathrm{NFU}^{-\mathrm{AC}}$ to denote $\mathrm{NFU}$ minus the Axiom of Choice. 

\begin{Definitions1}
We use $\mathcal{L}_{\mathrm{NFU}}$ to denote the extension of $\mathcal{L}$ obtained adding a unary predicate $\mathcal{S}$ and a binary function symbol $\langle \cdot, \cdot \rangle$.
\end{Definitions1}

The unary predicate $\mathcal{S}$ will be used to distinguish sets from urelements and $\langle \cdot, \cdot \rangle$ will act as a type-level pairing function. Before presenting the axioms of $\mathrm{NFU}^{-\mathrm{AC}}$ we first need to extend the notion of stratification to $\mathcal{L}_{\mathrm{NFU}}$-formulae.

\begin{Definitions1}
The terms of $\mathcal{L}_{\mathrm{NFU}}$ are built-up inductively from variables using the function $\langle \cdot, \cdot \rangle$. Let $\phi$ be an $\mathcal{L}_{\mathrm{NFU}}$-formula. We use $\mathbf{Term}(\phi)$ to denote the set of $\mathcal{L}_{\mathrm{NFU}}$-terms appearing in $\phi$. We say that $\sigma: \mathbf{Term}(\phi) \longrightarrow \mathbb{N}$ is a stratification of $\phi$ if for all terms $s$ and $t$ appearing in $\phi$, 
\begin{itemize}
\item[(i)] if $s$ is a term appearing in $t$ then $\sigma(\textrm{`}t\textrm{'})= \sigma(\textrm{`}s\textrm{'})$, 
\item[(ii)] if $s \in t$ is a subformula of $\phi$ then $\sigma(\textrm{`}t\textrm{'})= \sigma(\textrm{`}s\textrm{'})+1$,  
\item[(iii)] if $s = t$ is a subformula of $\phi$ then $\sigma(\textrm{`}t\textrm{'})= \sigma(\textrm{`}s\textrm{'})$. 
\end{itemize}
If there exists a stratification of $\phi$ then we say that $\phi$ is stratified. 
\end{Definitions1}

\begin{Definitions1}
We use $\mathrm{NFU}^{-\mathrm{AC}}$ to denote the $\mathcal{L}_{\mathrm{NFU}}$-theory with axioms:
\begin{itemize}
\item[](Weak Extensionality) $\forall x \forall y(\mathcal{S}(x) \land \mathcal{S}(y) \Rightarrow (x=y \iff \forall z(z \in x \iff z \in y)))$ 
\item[](Stratified Comprehension) for all stratified $\mathcal{L}_{\mathrm{NFU}}$-formulae $\phi(x, \vec{z})$,
$$\forall \vec{z} \exists y (\mathcal{S}(y) \land \forall x(x \in y \iff \phi(x, \vec{z})))$$ 
\item[](Pairing) $\forall x \forall y \forall z \forall w(\langle x, y \rangle= \langle w, z \rangle \Rightarrow (x=w \land y=z))$ 
\end{itemize}
\end{Definitions1}

\begin{Definitions1}
We use $\mathrm{NFU}$ to denote the $\mathcal{L}_{\mathrm{NFU}}$-theory that obtained from $\mathrm{NFU}^{-\mathrm{AC}}$ by adding the Axiom of Choice.
\end{Definitions1}

Following \cite{hol98} and \cite{hol01} we have opted to include an axiom that asserts the existence of a type-level pairing function in our axiomatisation of $\mathrm{NFU}^{-\mathrm{AC}}$. We will indicate below how this pairing function implies that there exists a Dedekind infinite set. Without the Axiom of Choice, Jensen's theory \cite{jen69} fortified with an axiom asserting the existence of a Dedekind infinite set is not equivalent to $\mathrm{NFU}^{-\mathrm{AC}}$, however they do have the same consistency strength. One way of seeing this is to use \cite{jen69} combined with work done in \cite{mat01} to see that the Axiom of Choice can be consistently added to Jensen's theory \cite{jen69} fortified with an axiom asserting the existence of a Dedekind infinite set. In the presence of both the Axiom of Choice and an axiom asserting the existence of a Dedekind infinite set, Jensen's theory \cite{jen69} is equipped with a type-level pairing function that can be used in instances of the comprehension scheme. This shows that $\mathrm{NFU}^{-\mathrm{AC}}$ and $\mathrm{TSTI}$ have exactly the same consistency strength:

\begin{Theorems1} \label{Th:ConsistencyStregthOfNFU}
(Jensen) $\mathrm{Con}(\mathrm{NFU}^{-\mathrm{AC}}) \iff \mathrm{Con}(\mathrm{TSTI})$. \Square
\end{Theorems1}

The set theory $\mathrm{NF}$ can be obtained from $\mathrm{NFU}^{-\mathrm{AC}}$ by adding an axiom which says that everything is a set.

\begin{Definitions1}
We use $\mathrm{NF}$ to denote the $\mathcal{L}_{\mathrm{NFU}}$-theory obtained from $\mathrm{NFU}^{-\mathrm{AC}}$ by adding the axiom 
\begin{equation} \label{Eq:EveryObjectIsASet}
\forall x(\mathcal{S}(x))
\end{equation} 
\end{Definitions1}

It should be noted that we could have axiomatised $\mathrm{NF}$, as is done in \cite{qui37}, in the language $\mathcal{L}$. In the presence of (\ref{Eq:EveryObjectIsASet}) the symbol $\mathcal{S}$ becomes redundant and the Weak Extensionality axiom reduces to the usual extensionality axiom for set theory. It follows from \cite{spe53} and \cite{qui45} that any model of $\mathrm{NF}$ can be expanded to a model with a pairing function $\langle \cdot, \cdot\rangle$ that satisfies the Pairing Axiom and can be used in instances of the Stratified Comprehension scheme without raising types. It should also be noted that \cite{spe53} shows that the Axiom of Choice is inconsistent with $\mathrm{NFU}^{-\mathrm{AC}}$ plus (\ref{Eq:EveryObjectIsASet}).\\
\\
\indent Stratified Comprehension in the theory $\mathrm{NFU}^{-\mathrm{AC}}$ guarantees the existence of a universal set which we denote $V$. The fact that the function $x \mapsto \langle x, x \rangle$ is injective but not surjective implies that $V$ is Dedekind infinite. Cardinal and ordinal numbers in $\mathrm{NFU}^{-\mathrm{AC}}$ are defined to be equivalence classes of equipollent sets and equivalence classes of isomorphic well-orderings respectively. If $X$ is a set then we use $|X|$ to denote the cardinal number such that $X \in |X|$. Stratified Comprehension ensures that both the set of all ordinals ($\mathrm{NO}$) and the set of all cardinals ($\mathrm{NC}$) exist. We use $\mathrm{NCI}$ to denote the set of infinite cardinals. The least cardinal number, denoted $\mathbf{0}$, is the set of all sets and urelements that have no members. We use $\iota$ denote the function $x \mapsto \{x\}$. Equipped with a successor operation ($S$) we are able to define the natural numbers ($\mathbb{N}$) as the smallest inductive set:  
$$S(x)= \{y \mid (\exists z \in y)(y \backslash \iota z \in x)\}$$
$$\mathbb{N}= \bigcap \{x \mid (\mathbf{0} \in x) \land (\forall y \in x)(S(y) \in x)\}.$$
Define $+$ and $\cdot$ on $\mathbb{N}$ by: for all $k, m, n \in \mathbb{N}$,
$$k+m= n \textrm{ if there exists } x \in k \textrm { and } y \in m \textrm{ with } x \cap y = \emptyset \textrm{ and } n= |x \cup y|$$
$$k \cdot m= n \textrm{ if there exists } x \in k \textrm{ and } y \in m \textrm{ such that } |x \times y|= n.$$
By letting $\mathbf{1}= S(\mathbf{0})$ we obtain an $\mathcal{L}_\mathrm{PA}$-structure $\langle \mathbb{N}, +, \cdot, \mathbf{0}, \mathbf{1} \rangle$ that is a model of $\mathrm{PA}$. If $n$ is a concrete natural number then by adjoining the sets $\mathcal{P}(\mathbb{N}), \mathcal{P}^2(\mathbb{N}), \ldots, \mathcal{P}^{n-1}(\mathbb{N})$ to the structure $\langle \mathbb{N}, +, \cdot, \mathbf{0}, \mathbf{1} \rangle$ we obtain a structure that is a model of $n^{\textrm{th}}$ order arithmetic.\\
\\
\indent This interpretation of Peano Arithmetic allows $\mathrm{NFU}^{-\mathrm{AC}}$ to describe the syntax of recursive languages. If $\mathcal{L}^\prime$ is a recursive language then expressions in $\mathcal{L}^\prime$ can be coded as elements of $\mathbb{N}$, called a G\"{o}del coding, in such a way so as effective properties of $\mathcal{L}^\prime$ expressions are definable by arithmetic (and therefore stratified) $\mathcal{L}_{\mathrm{NFU}}$-formulae. Given a recursive language $\mathcal{L}^\prime$ we assume that a G\"{o}del coding of $\mathcal{L}^\prime$ has been fixed and we write $\ulcorner \phi \urcorner$ for the G\"{o}del code of $\phi$. We will often omit the corners and equate a formula $\phi$ with its G\"{o}del code. In \cite{hin75} Hinnion shows that if an $\mathcal{L}$-structure $\mathcal{M}$ is a set then there is a stratified formula $\mathrm{Sat}_\mathcal{L}(\phi, a, \mathcal{M})$ which says that $\phi$ is an $\mathcal{L}$-formulae, $a$ is sequence of elements of $M$ that agrees with the arity of $\phi$ and $\mathcal{M}$ satisfies $\phi[a]$. If $\mathcal{L}^\prime$ is a recursive language then Hinnion's definition of satisfaction for $\mathcal{L}$-structures can easily be extended to define a ternary stratified formula $\mathrm{Sat}_{\mathcal{L}^\prime}$ which expresses satisfaction in an $\mathcal{L}^\prime$-structure. Using the stratified formulae $\mathrm{Sat}_{\mathcal{L}_\mathrm{PA}}$ one can see that $\mathrm{NFU}^{-\mathrm{AC}}$ proves the single sentence which asserts that the structure $\langle \mathbb{N}, +, \cdot, \mathbf{0}, \mathbf{1} \rangle$ is a model of $\mathrm{PA}$. And, moreover, for any concrete natural number $n$, $\mathrm{NFU}^{-\mathrm{AC}}$ proves the single sentence which asserts that the structure $\langle \mathcal{P}^{n-1}(\mathbb{N}), \ldots, \mathcal{P}(\mathbb{N}), \mathbb{N}, +, \cdot, \mathbf{0}, \mathbf{1} \rangle$ is a model of $n^\textrm{th}$ order arithmetic. If $\mathcal{M}$ is a set structure in a recursive language $\mathcal{L}^\prime$, $\phi$ is an $\mathcal{L}^\prime$-formula and $a$ is a sequence of elements of $\mathcal{M}$ then we will write $\mathcal{M} \models \phi[a]$ instead of $\mathrm{Sat}_{\mathcal{L}^\prime}(\phi, a, \mathcal{M})$.\\
\\             
\indent The following definition mirrors the definition of an initial ordinal in $\mathrm{ZFC}$:

\begin{Definitions1}
We say that an ordinal $\alpha$ is initial if 
$$(\forall \delta < \alpha)(|\{\beta \mid \beta < \delta\}| < |\{\beta \mid \beta < \alpha\}|).$$
We use $\omega$ to denote the first initial ordinal, $\omega_1$ to denote the least initial ordinal $>\omega$, $\omega_2$ to denote the least initial ordinal $>\omega_1$, etc.
\end{Definitions1}

In $\mathrm{ZFC}$ cardinals correspond to initial ordinals. It is important to note that in $\mathrm{NFU}^{-\mathrm{AC}}$ this coincidence does not occur. If $R$ is a binary relation then we will write $\mathrm{Dom}(R)$ for $\mathrm{dom}(R) \cup \mathrm{rng}(R)$.

\begin{Definitions1}
Let $\alpha$ be an ordinal. Define
$$\mathrm{Card}(\alpha)= |\mathrm{Dom}(R)| \textrm{ where } R \in \alpha.$$
For all $n \in \mathbb{N}$, define $\aleph_n= \mathrm{Card}(\omega_n)$.
\end{Definitions1} 

\indent One unorthodox feature of $\mathrm{NFU}^{-\mathrm{AC}}$ is the fact that it proves that there are sets, for example $V$, which are not the same size as their own set of singletons. This motivates the introduction of the $T$ operation which is defined on cardinals, ordinals and equivalence classes of isomorphic well-founded relations, and the definition of Cantorian and strongly Cantorian sets. If $R$ is a relation then we use $[R]$ to denote the set of all relations isomorphic to $R$. If $F$ is a function and $X$ is a set then we write $F``X$ for the set of all things that can be obtained by applying $F$ to an element of $X$.

\begin{Definitions1}
We say that a set $X$ is Cantorian if $|X|= |\iota``X|$. We say a set $X$ is strongly Cantorian if the restriction of the map $\iota$ to $X$ witnesses the fact that $|X|= |\iota``X|$. 
\end{Definitions1}

\begin{Definitions1}
If $X$ is a set then define
$$T(|X|)= |\iota``X|.$$
If $R$ is a well-founded relation then define
$$T([R])= [\{\langle \iota x, \iota y \rangle \mid \langle x, y \rangle \in R\}].$$ 
\end{Definitions1}

The $T$ operation commutes with the functions $+$ and $\cdot$ defined on $\mathbb{N}$ and is the identity on $\mathbf{0}$ and $\mathbf{1}$. Since both $T``\mathbb{N}$ and $T^{-1}``\mathbb{N}$ contain $\mathbf{0}$ and are closed under $S$, it follows that $T``\mathbb{N}= \mathbb{N}$. Therefore the $T$ operation is an automorphism of the interpretation of arithmetic in a model of $\mathrm{NFU}^{-\mathrm{AC}}$. The Axiom of Counting ($\mathrm{AxCount}$) asserts that this automorphism is the identity:
\begin{itemize}
\item[]($\mathrm{AxCount}$) $(\forall n \in \mathbb{N})(T(n)=n)$
\end{itemize}
This axiom was first introduced by J. Barkley Rosser in \cite{ros78} in order to facilitate induction in $\mathrm{NF}$. Steve Orey's \cite{ore64} shows that $\mathrm{NF}+\mathrm{AxCount}$ proves $\mathrm{Con}(\mathrm{NF})$. As part of \cite{hol01}, which also initiates the comparison of extensions of $\mathrm{NFU}$ with subsystems and extensions of $\mathrm{ZFC}$, Randall Holmes investigates the strength of the theory $\mathrm{NFU}+\mathrm{AxCount}$ in terms of which infinite cardinals this theory proves exist. In \cite{for77} Thomas Forster identifies two natural weakenings of $\mathrm{AxCount}$:
\begin{itemize}
\item[]($\mathrm{AxCount}_\leq$) $(\forall n \in \mathbb{N})(n \leq T(n))$ 
\item[]($\mathrm{AxCount}_\geq$) $(\forall n \in \mathbb{N})(n \geq T(n))$
\end{itemize} 
Many of the strong consequences of $\mathrm{AxCount}$ also follow from the weaker assumption $\mathrm{AxCount}_\leq$. For example, \cite{hin75} shows that $\mathrm{NF}+\mathrm{AxCount}_\leq$ proves the consistency of Zermelo Set Theory. And \cite{fh09} shows that if $\mathrm{NFU}^{-\mathrm{AC}}+\mathrm{AxCount}_\leq$ is consistent then so is $\mathrm{NFU}^{-\mathrm{AC}}+ (\textrm{the function on } \mathbb{N} \textrm{ defined by } n \mapsto V_n \textrm{ exists})$. In contrast it is not known if $\mathrm{NF}+\mathrm{AxCount}_\geq$ proves the consistency of Zermelo Set Theory. And the assertion that the function on $\mathbb{N}$ defined by $n \mapsto V_n$ exists proves $\mathrm{AxCount}_\leq$. The relative strengths of $\mathrm{AxCount}$, $\mathrm{AxCount}_\leq$ and $\mathrm{AxCount}_\geq$ over $\mathrm{NFU}$ is studied in \cite{mck15}:

\begin{Theorems1} \label{Th:VariantsOfAxCountOverNFU}
\begin{itemize}
\item[(I)] $\mathrm{NFU}+\mathrm{AxCount} \vdash \mathrm{Con}(\mathrm{NFU}+\mathrm{AxCount}_\leq)$
\item[(II)] $\mathrm{NFU}+\mathrm{AxCount}_\leq \vdash \mathrm{Con}(\mathrm{NFU}+\mathrm{AxCount}_\geq)$
\end{itemize}
\end{Theorems1}

\cite{mck15} also provides evidence which appears to indicate that, over $\mathrm{NFU}$, $\mathrm{AxCount}_\geq$ is weak.

\begin{Theorems1} \label{Th:EvidenceThatAxCountGEQIsWeak}
There is a model of $\mathrm{NFU}+\mathrm{AxCount}_\geq$ which believes that $\mathrm{NCI}$ is finite.
\end{Theorems1}

\cite{mat01} shows that $\mathrm{TSTI}$ is equiconsistent with the set theory $\mathrm{MOST}$ that is a subsystem of $\mathrm{ZFC}$ that includes the Axiom of Choice. The fact that equiconsistencies between $\mathrm{MOST}$ and $\mathrm{TSTI}$, and $\mathrm{TSTI}$ and $\mathrm{NFU^{-\mathrm{AC}}}$ show that these theories have the same arithmetic yields the following strong equiconsistency:

\begin{Theorems1} \label{Th:StrongConsistencyOfChoiceWithNFU}
(Jensen, Mathias) 
$$\mathrm{Con}(\mathrm{NFU}^{-\mathrm{AC}}) \iff \mathrm{Con}(\mathrm{NFU}).$$
Moreover, if $\phi$ is an $\mathcal{L}_\mathrm{PA}$-sentence then
$$\mathrm{NFU}^{-\mathrm{AC}} \vdash \phi \textrm{ if and only if } \mathrm{NFU} \vdash \phi.$$ 
\end{Theorems1} 

It follows from Theorem \ref{Th:StrongConsistencyOfChoiceWithNFU} that any occurrence of ``$\mathrm{NFU}$" in Theorem \ref{Th:VariantsOfAxCountOverNFU} can be replaced by ``$\mathrm{NFU}^{-\mathrm{AC}}$".\\
\\
\indent In \cite{hin75} Hinnion shows that subsystems of $\mathrm{ZFC}$ can be interpreted in substructures of the set of equivalence classes of isomorphic topped well-founded extensional relations in $\mathrm{NF}$. Hinnion's techniques have since been established (see \cite{hol98}, \cite{hol01} and \cite{solXX}) as the standard method for proving lower bounds on the consistency strength of extensions of $\mathrm{NFU}^{-\mathrm{AC}}$ relative to subsystems and extensions of $\mathrm{ZFC}$.  

\begin{Definitions1} \label{Df:BFEXT}
A structure $\langle A, R \rangle$, where $R$ is a binary relation, is a $\mathrm{BFEXT}$ if
\begin{itemize}
\item[(i)] $\mathrm{Dom}(R)= A$, 
\item[(ii)] $\forall B((B \neq \emptyset \land B \subseteq A) \Rightarrow (\exists b \in B)(\forall c \in B) \neg (\langle c, b \rangle \in R))$,
\item[(iii)] $(\forall a, b \in A)(a=b \iff \forall c(\langle c, a \rangle \in R \iff \langle c, b \rangle \in R))$.
\end{itemize}
We say that a binary relation $R$ is a $\mathrm{BFEXT}$ if $\langle \mathrm{Dom}(R), R \rangle$ is a $\mathrm{BFEXT}$. 
\end{Definitions1}

\begin{Definitions1} \label{Df:Segment}
Let $\langle A, R \rangle$ be a $\mathrm{BFEXT}$. If $a \in A$ then define
$$\mathrm{seg}_R(a)= R \upharpoonright \bigcap \{B \subseteq A \mid (a \in B \land (\forall b \in B)(\forall c \in B)(\langle c, b \rangle \in R \Rightarrow c \in B))\}.$$
\end{Definitions1}

\begin{Definitions1}
$$\Omega= \{R \mid (R \textrm{ is a } \mathrm{BFEXT}) \land (\exists a \in \mathrm{Dom}(R))(R= \mathrm{seg}_R(a))\}$$
\end{Definitions1}

The fact that $\Omega$ is defined by a stratified set abstract shows that $\mathrm{NFU}^{-\mathrm{AC}}$ proves that $\Omega$ is a set. If $R \in \Omega$ then the $a \in \mathrm{Dom}(R)$ with $R= \mathrm{seg}_R(a)$ is unique--- we will use $\mathds{1}_R$ to denote this element. We will sometimes call $[R]$ the type of $R$.

\begin{Definitions1}
The structure $\langle \mathrm{BF}, \mathcal{E} \rangle$ is defined by 
$$\mathrm{BF}= \{[R] \mid R \in \Omega\}$$
$$\mathcal{E}= \{\langle [R], [S] \rangle \in \mathrm{BF} \times \mathrm{BF} \mid (\exists a \in \mathrm{Dom}(S))(R \cong \mathrm{seg}_S(a) \land \langle a, \mathds{1}_S \rangle \in S)\}.$$
\end{Definitions1}

\begin{Theorems1} \label{Th:BFisWellFoundedAndExtensional}
(Hinnion) The structure $\langle \mathrm{BF}, \mathcal{E} \rangle$ is well-founded and extensional. \Square
\end{Theorems1}

A consequence of this theorem is that if $\mathfrak{a}$ is an equivalence class of isomorphic BFEXTs then $\mathrm{seg}_{\mathcal{E}}(\mathfrak{a})$ is a BFEXT. The BFEXTs represented by $\mathfrak{a}$ are related to $\mathrm{seg}_{\mathcal{E}}(\mathfrak{a})$ by the $T$ operation.

\begin{Lemma1} \label{Th:RelationshipBetweenBFEXTAndSeg}
(Hinnion) If $\mathfrak{a} \in \mathrm{BF}$ then $[\mathrm{seg}_{\mathcal{E}}(\mathfrak{a})]= T^2(\mathfrak{a})$. \Square
\end{Lemma1}

By considering rank initial segments of the structure $\langle \mathrm{BF}, \mathcal{E} \rangle$ Hinnion builds models of subsystems of $\mathrm{ZFC}$.

\begin{Definitions1}
Let $S$ be a well-founded extensional relation. Define
$$S^0= \{a \in \mathrm{Dom}(S) \mid \neg(\exists b \in \mathrm{Dom}(S))(\langle b, a \rangle \in S)\}$$
$$S^{\alpha+1}= \{a \in \mathrm{Dom}(S) \mid \forall b(\langle b, a\rangle \in S \Rightarrow b \in S^\alpha)\}$$
$$S^\lambda= \bigcup_{\alpha<\lambda} S^\alpha \textrm{ for limit } \lambda.$$
\end{Definitions1}  

\noindent Note that for well-founded extensional $S$, the formula `$x=S^\alpha$' is stratified and admits a stratification which assigns the same type to the variables `$x$' and `$\alpha$'.

\begin{Definitions1}
For an ordinal $\alpha$, we use $M_\alpha$ to denote $\mathcal{E}^{\omega+\alpha}$.
\end{Definitions1}

\noindent Note that $M_0= \{[R] \mid |\mathrm{Dom}(R)|< \aleph_0\}$ and $|M_0|= \aleph_0$.\\
\\    
\indent One of the achievements of \cite{hin75}, which we mentioned above, was to show that if $\mathrm{AxCount}_\leq$ holds in $\mathrm{NF}$ then the single sentence asserting that $M_{\omega}$ is a model of Zermelo Set Theory is provable. Even though the setting of \cite{hin75} is $\mathrm{NF}$, Hinnion's argument can also be carried out in the weaker theory $\mathrm{NFU}^{-\mathrm{AC}}$ (see \cite{hol98} and \cite{solXX}). Combining this with the work in \cite{mat01} on the consistency strength of $\mathrm{TSTI}$ we note the following weak version of Hinnion's result which we will use in the next section.

\begin{Theorems1} \label{Th:HinnionLowerBoundOnStrengthOfAxCountLEQ}
(Hinnion) $\mathrm{NFU}^{-\mathrm{AC}} + \mathrm{AxCount}_\leq \vdash \mathrm{Con}(\mathrm{TSTI})$. \Square
\end{Theorems1}

\section[$\mathrm{NFU}^{-\mathrm{AC}}+\mathrm{AxCount}_\geq$ proves the consistency of $\mathrm{TSTI}$]{$\mathrm{NFU}^{-\mathrm{AC}}+\mathrm{AxCount}_\geq$ proves the consistency of $\mathrm{TSTI}$}

In this section I will show that $\mathrm{NFU}^{-\mathrm{AC}}+\mathrm{AxCount}_\geq$ proves the consistency of $\mathrm{TSTI}$. The main tool used to prove this result is the technique, developed in \cite{hin75}, of using the class of topped well-founded extensional relations in $\mathrm{NFU}^{-\mathrm{AC}}$ to interpret well-founded set theories. It follows from Theorem \ref{Th:HinnionLowerBoundOnStrengthOfAxCountLEQ} that if $\mathrm{AxCount}$ holds in $\mathrm{NFU}^{-\mathrm{AC}}$ then $\mathrm{Con}(\mathrm{TSTI})$ holds. In light of this, all we need to prove is that the theory $\mathrm{NFU}^{-\mathrm{AC}}+\mathrm{AxCount}_\geq+\neg \mathrm{AxCount}$ proves the consistency of $\mathrm{TSTI}$.\\
\\
\indent Let $\mathcal{M}$ be a model of $\mathrm{NFU}^{-\mathrm{AC}}+\mathrm{AxCount}_\geq+ \neg\mathrm{AxCount}$. The proof will show that there is an elementary $\mathcal{L}_{\mathrm{PA}}$-substructure $\mathcal{A}$ of $\langle \mathbb{N}^\mathcal{M}, +^\mathcal{M}, \cdot^\mathcal{M}, \mathbf{0}^\mathcal{M}, \mathbf{1}^\mathcal{M} \rangle$ that satisfies $\mathrm{Con}(\mathrm{TSTI})$. It then follows from the elementarity of $\mathcal{A}$ that $\mathcal{M}$ satisfies $\mathrm{Con}(\mathrm{TSTI})$. The fact that $\mathcal{A}$ satisfies $\mathrm{Con}(\mathrm{TSTI})$ will be obtained by showing that for every $n \in \mathcal{A}$, there is a set substructure of $\langle \mathrm{BF}, \mathcal{E} \rangle$ that is a model of $\mathrm{TSTI}_n$. The structure $\mathcal{A}$ is obtained by considering the fixed points of $T$ acting on $\mathbb{N}^\mathcal{M}$. The following Lemma is proved in the  theory $\mathrm{NFU}^{-\mathrm{AC}}+\mathrm{AxCount}_\geq+ \neg\mathrm{AxCount}$:

\begin{Lemma1} \label{Th:FixedPointSetOfTIsInitialSegment}
If $n \in \mathbb{N}$ is such that $T(n)=n$ then for all $m \leq n$, $T(m)=m$.
\end{Lemma1}

\begin{proof}
Suppose that there are $n, m \in \mathbb{N}$ with $m < n$, $T(n)=n$ and $T(m) < m$. But then $n-m \in \mathbb{N}$ and
$$T(n-m)= T(n) - T(m) = n - T(m) > n - m$$
which contradicts $\mathrm{AxCount}_\geq$.
\Square
\end{proof}

\begin{Definitions1}
Define $\mathcal{A}= \langle A, +^\mathcal{A}, \cdot^\mathcal{A}, \mathbf{0}^\mathcal{A}, \mathbf{1}^\mathcal{A} \rangle$ to be the $\mathcal{L}_\mathrm{PA}$-substructure of\\ 
$\langle \mathbb{N}^\mathcal{M}, +^\mathcal{M}, \cdot^\mathcal{M}, \mathbf{0}^\mathcal{M}, \mathbf{1}^\mathcal{M} \rangle$ with domain
$$A=\{n \in \mathbb{N}^\mathcal{M} \mid \mathcal{M} \models (T(n)=n)\}.$$
\end{Definitions1}

Note that since $A \subseteq \mathbb{N}^\mathcal{M}$, $A \neq \mathbb{N}$, $\mathbf{0}^\mathcal{M} \in A$ and $A$ is closed under $S^\mathcal{M}$ it follows that $A$ is not a set of $\mathcal{M}$.

\begin{Lemma1}
The $\mathcal{L}_\mathrm{PA}$-structure $\langle \mathbb{N}^\mathcal{M}, +^\mathcal{M}, \cdot^\mathcal{M}, \mathbf{0}^\mathcal{M}, \mathbf{1}^\mathcal{M} \rangle$ is a proper elementary end-extension of $\mathcal{A}$. 
\end{Lemma1}

\begin{proof}
Lemma \ref{Th:FixedPointSetOfTIsInitialSegment} implies that  $\langle \mathbb{N}^\mathcal{M}, +^\mathcal{M}, \cdot^\mathcal{M}, \mathbf{0}^\mathcal{M}, \mathbf{1}^\mathcal{M} \rangle$ is an end-extension of $\mathcal{A}$. It follows from the fact that $\mathcal{M} \models \neg \mathrm{AxCount}$ that $A \neq \mathbb{N}^\mathcal{M}$. That $\mathcal{A} \prec \langle \mathbb{N}^\mathcal{M}, +^\mathcal{M}, \cdot^\mathcal{M}, \mathbf{0}^\mathcal{M}, \mathbf{1}^\mathcal{M} \rangle$ follows since $T^\mathcal{M}$ is an automorphism of the structure  $\langle \mathbb{N}^\mathcal{M}, +^\mathcal{M}, \cdot^\mathcal{M}, \mathbf{0}^\mathcal{M}, \mathbf{1}^\mathcal{M} \rangle$ and $\mathrm{PA}$ has definable Skolem functions. 
\Square
\end{proof}

We now turn to showing that for all $n \in A$, $\mathcal{A}$ satisfies $\mathrm{Con}(\mathrm{TSTI}_n)$. This will be achieved by working in $\mathcal{M}$ and showing that there is an $n \in \mathbb{N}$ such that $T(n) < n$ and there exists a set model of $\mathrm{TSTI}_n$ in the structure $\langle \mathrm{BF}, \mathcal{E} \rangle$. From this point on we work inside $\mathcal{M}$. 

\begin{Definitions1}
If $\kappa$ is a cardinal then define
$$2^{\kappa}= \left\{\begin{array}{ll}
|\mathcal{P}(X)| & \textrm{if there exists a set } X \textrm{ with } |\iota``X|= \kappa\\
\emptyset & \textrm{otherwise}
\end{array} \right)$$
\end{Definitions1}

\noindent This modification of the usual definition of cardinal exponentiation has the property that the function $\kappa \mapsto 2^\kappa$ is definable by a stratified formula which admits a stratification that assigns the same type to the result and the argument of the function. The following result shows that this exponentiation operation possesses the strictly inflationary property that we intuitively associate with cardinal exponentiation.

\begin{Lemma1} \label{Th:ExpInflationary}
Let $\kappa$ be a cardinal. If $2^{\kappa} \neq \emptyset$ then $\kappa < 2^{\kappa}$.
\end{Lemma1}

\begin{proof}
The usual proof of Cantor's Theorem yields for all $X$, $|\iota``X| < |\mathcal{P}(X)|$ and this proof only appeals stratified instances of comprehension. 
\Square
\end{proof}

\begin{Definitions1}
Define $\beth: \mathbb{N} \longrightarrow \mathrm{NC}$ by
$$\beth(0)= \aleph_0$$
$$\beth(n+1)= \left\{\begin{array}{ll}
2^{\beth(n)} & \textrm{if } \beth(n) \in \mathrm{NC}\\
\emptyset & \textrm{if } \beth(n)= \emptyset
\end{array}\right.$$
\end{Definitions1}

\noindent Note that stratified comprehension ensures that $\beth$ is a set. The following results are proved or adapted from results proved in \cite{hin75}.

\begin{Lemma1}
If $2^\kappa \neq \emptyset$ then $T(2^\kappa)=2^{T(\kappa)}$.
\end{Lemma1}

\begin{proof}
This follows immediately from the fact that for all $X$, $|\iota``\mathcal{P}(X)|= |\mathcal{P}(\iota``X)|$. 
\Square
\end{proof}

\begin{Lemma1} \label{Th:TFunctionOnBeths}
Let $n \in \mathbb{N}$. If $\beth(n) \neq \emptyset$ then $\beth(T(n)) = T(\beth(n))$.
\end{Lemma1}

\begin{proof}
We prove this by induction on $n$. It holds for $n=0$. Suppose that the Lemma holds for $n$. Suppose that $\beth(n+1) \neq \emptyset$. Therefore $\beth(n) \neq \emptyset$ and so $\beth(T(n))= T(\beth(n))$. So,
$$\beth(T(n+1))= \beth(T(n)+1)= 2^{\beth(T(n))}= 2^{T(\beth(n))}= T(2^{\beth(n)})= T(\beth(n+1)).$$
\Square 
\end{proof}

\begin{Lemma1} \label{Th:ConditionForPowerOf2ToExist}
If $\kappa \leq T(|V|)$ then $2^\kappa \neq \emptyset$.
\end{Lemma1}

\begin{proof}
Let $X \in \kappa$ and let $f: X \longrightarrow \iota``V$ be an injection. Let $B= \mathrm{rng}(f)$ and let $A= \bigcup B$. Therefore $|\iota``A|= \kappa$. \Square
\end{proof}

\begin{Lemma1} \label{Th:ExistenceOfNonStandardN}
There exists an $n \in \mathbb{N}$ with $T(n) < n$ such that $\beth(n) \neq \emptyset$ and $\beth(n) \leq T^4(|V|)$.
\end{Lemma1}

\begin{proof}
If $\beth(n) \neq \emptyset$ and $\beth(n) \leq T^4(|V|)$ for all $n \in \mathbb{N}$, then we are done since $\mathrm{AxCount}$ fails and $\mathrm{AxCount}_\geq$ holds. Suppose that $n+1$ is least such that $\beth(n+1) = \emptyset$ or $\beth(n+1) \nleq T^4(|V|)$. If $T(n+1)= n+1$ then $T(n)= n$. And, $\beth(n) \neq \emptyset$ and $\beth(n) \leq T^4(|V|)$. But then, by Lemma \ref{Th:ConditionForPowerOf2ToExist}, $\beth(n+1) \neq \emptyset$. And, by Lemma \ref{Th:TFunctionOnBeths}, $T(\beth(n+1))= \beth(n+1)$. So, $\beth(n+1) \leq |V|$ implies that $\beth(n+1) \leq T^4(|V|)$, which contradicts our assumptions. Therefore $T(n+1) < n+1$ and $T(n) < n$. Since $n+1$ was least such that $\beth(n+1) = \emptyset$ or $\beth(n+1) \nleq T^4(|V|)$, it follows that $\beth(n) \neq \emptyset$ and $\beth(n) \leq T^4(|V|)$.      
\Square
\end{proof}

\begin{Lemma1} \label{Th:SizeOfM}
Let $n \in \mathbb{N}$. If $\beth(T(n)) \neq \emptyset$ then $|M_n| \leq \beth(T(n))$. 
\end{Lemma1}

\begin{proof}
The assertion `$\beth(T(n)) \neq \emptyset \Rightarrow |M_n| \leq \beth(T(n))$' is stratified, so we can prove it by induction. The base case holds because $|M_0|= \aleph_0$. Suppose that the Lemma holds for some $n \in \mathbb{N}$ and assume that $\beth(T(n+1)) \neq \emptyset$. Note that if $\mathfrak{a} \in M_{n+1}$ then $\{\mathfrak{b} \mid \langle \mathfrak{b}, \mathfrak{a} \rangle \in \mathcal{E}\} \subseteq M_n$. The map $g: \iota``M_{n+1} \longrightarrow \mathcal{P}(M_n)$ defined by $\{\mathfrak{a}\} \mapsto \{\mathfrak{b} \mid \langle \mathfrak{b}, \mathfrak{a} \rangle \in \mathcal{E}\}$ is injective. Therefore $T(|M_{n+1}|) \leq |\mathcal{P}(M_n)|$. Now, $|\iota``M_n|= T(|M_n|)$, therefore $2^{T(|M_n|)}= |\mathcal{P}(M_n)|$. We also know that $2^{\beth(T(n))} \neq \emptyset$. And so,
$$T(|M_{n+1}|) \leq 2^{T(|M_n|)} \leq 2^{T(\beth(T(n)))}= T(2^{\beth(T(n))})= T(\beth(T(n+1))).$$  
\Square
\end{proof}

\begin{Lemma1} \label{Th:EveryBFEXTisT2}
Let $n \in \mathbb{N}$ be such that $\beth(n) \neq \emptyset$ and $\beth(n) \leq T^4(|V|)$. If $\mathfrak{a} \in M_n$ then there is a $\mathfrak{b} \in \mathrm{BF}$ such that $\mathfrak{a}= T^2(\mathfrak{b})$. 
\end{Lemma1}

\begin{proof}
Let $n \in \mathbb{N}$ be such that $\beth(n) \neq \emptyset$ and $\beth(n) \leq T^4(|V|)$. Let $\mathfrak{a} \in M_n$. Therefore $\mathrm{Dom}(\mathrm{seg}_\mathcal{E}(\mathfrak{a})) \subseteq M_n$, and so $|\mathrm{Dom}(\mathrm{seg}_\mathcal{E}(\mathfrak{a}))| \leq |M_n|$. Since $T(n) \leq n$, we know that $\beth(T(n)) \neq \emptyset$. Therefore, by Lemma \ref{Th:SizeOfM},
$$|\mathrm{Dom}(\mathrm{seg}_\mathcal{E}(\mathfrak{a}))| \leq \beth(T(n)) \leq \beth(n) \leq T^4(|V|).$$
Let $f: \mathrm{Dom}(\mathrm{seg}_\mathcal{E}(\mathfrak{a})) \longrightarrow \iota^4``V$ be an injection. Let $A= \bigcup^4 \mathrm{rng}(f)$ and define $S \subseteq A \times A$ by
$$\langle x, y \rangle \in S \textrm{ if and only if } \langle f^{-1}(\iota^4 x), f^{-1}(\iota^4 y) \rangle \in \mathrm{seg}_\mathcal{E}(\mathfrak{a}).$$
Let $\mathfrak{b}= [S]$. By Lemma \ref{Th:RelationshipBetweenBFEXTAndSeg} we have
$$T^2(\mathfrak{a})= [\mathrm{seg}_\mathcal{E}(\mathfrak{a})]= T^4(\mathfrak{b}).$$
Therefore $\mathfrak{a}= T^2(\mathfrak{b})$.
\Square
\end{proof}

\begin{Lemma1} \label{Th:KeyCompLemma}
Let $n \in \mathbb{N}$ be such that $\beth(n) \neq \emptyset$ and $\beth(n) \leq T^4(|V|)$. If $B \subseteq M_n$ then there exists $\mathfrak{d} \in M_{n+1}$ such that for all $\mathfrak{b} \in M_n$,
$$\langle \mathfrak{b}, \mathfrak{d} \rangle \in \mathcal{E} \textrm{ if and only if } \mathfrak{b} \in B.$$ 
\end{Lemma1}

\begin{proof}
Let $B \subseteq M_n$. It follows from Lemma \ref{Th:EveryBFEXTisT2} that for all $\mathfrak{b}$, there exists $\mathfrak{b}^\prime \in \mathrm{BF}$ such that $\mathfrak{b}= T^2(\mathfrak{b}^\prime)$. Therefore, by Lemma \ref{Th:RelationshipBetweenBFEXTAndSeg}, for all $\mathfrak{b} \in \mathrm{BF}$, $\mathfrak{b}= [\mathrm{seg}_{\mathcal{E}}(\mathfrak{b}^\prime)]$. Define
$$S= \left(\bigcup \{ \mathrm{seg}_{\mathcal{E}}(\mathfrak{b}^\prime) \mid T^2(\mathfrak{b}^\prime) \in B\}\right) \cup \{\langle \mathfrak{b}^\prime, V \rangle \mid T^2(\mathfrak{b}^\prime) \in B\}.$$
Now, $S \in \Omega$ and so let $\mathfrak{d}= [S]$. It is clear that $\mathfrak{d} \in M_{n+1}$ and for all $\mathfrak{b} \in M_n$,
$$\langle \mathfrak{b}, \mathfrak{d} \rangle \in \mathcal{E} \textrm{ if and only if } \mathfrak{b} \in B.$$
\Square 
\end{proof}

Equipped with these results we are now in a position to show that $\mathcal{M}$ satisfies $\mathrm{Con}(\mathrm{TSTI})$.

\begin{Lemma1}
$\mathcal{M} \models \mathrm{Con}(\mathrm{TSTI})$.
\end{Lemma1}

\begin{proof}
Suppose that $\mathcal{M} \models \neg \mathrm{Con}(\mathrm{TSTI})$. Therefore
$$\mathcal{M} \models (\exists k \geq 4) \neg \mathrm{Con}(\mathrm{TSTI}_k).$$
Since the arithmetic of $\mathcal{M}$ is elementarily equivalent to $\mathcal{A}$ this implies
$$\mathcal{A} \models (\exists k \geq 4) \neg \mathrm{Con}(\mathrm{TSTI}_k).$$
Let $k \in A$ be such that $\mathcal{A} \models \neg \mathrm{Con}(\mathrm{TSTI}_k)$. Since $\mathcal{A}$ is an elementary submodel of the arithmetic of $\mathcal{M}$, this means that $\mathcal{M} \models \neg \mathrm{Con}(\mathrm{TSTI}_k)$. Note that $\mathcal{M} \models (T(k)=k)$.\\
Now, work inside $\mathcal{M}$. Let $n \in \mathbb{N}$ be such that $T(n) < n$, $\beth(n) \neq \emptyset$ and $\beth(n) \leq T^4(|V|)$. Lemma \ref{Th:ExistenceOfNonStandardN} ensures that there exists an $n \in \mathbb{N}$ with these properties. Note that for all $i \leq n$, $\beth(i) \neq \emptyset$ and $\beth(i) \leq T^4(|V|)$. We will build a set model of $\mathrm{TSTI}_{n+1}$. Since $k < n$ this will yield a contradiction. Let $D: [n] \longrightarrow V$ be defined by $D(i)= M_i$ for all $i \leq n$. Let $\in^\mathcal{D}: (n) \longrightarrow V$ be defined by $\in^\mathcal{D}(i)= \mathcal{E} \upharpoonright M_i \times M_{i+1}$. Note that Stratified Comprehension ensures that the functions $D$ and $\in^\mathcal{D}$ are sets. The structure $\mathcal{D}= \langle D, \in^\mathcal{D} \rangle$ is an $\mathcal{L}_{n+1}$-structure.\\  
We need to show that $\mathcal{D} \models \mathrm{TSTI}_{n+1}$. To see that $\mathcal{D} \models (\textrm{Extensionality})$ observe that the structure $\langle \mathrm{BF}, \mathcal{E} \rangle$ is extensional (Theorem \ref{Th:BFisWellFoundedAndExtensional}) and for each $i \leq n$, $\langle \mathrm{BF}, \mathcal{E} \rangle$ is an end-extension of $\langle M_i, \mathcal{E} \rangle$.\\
We now turn to showing that $\mathcal{D} \models (\textrm{Comprehension})$. Let $i < n$. Let $\phi(x^i, z_0^{j_0}, \ldots, z_{m-1}^{j_{m-1}})$ be an $\mathcal{L}_{n+1}$-formula (according to $\mathcal{M}$). We need to show that 
$$\mathcal{D} \models \forall \vec{z} \exists y^{i+1} \forall x^i(x^i \in_i y^{i+1} \iff \phi(x^i, \vec{z})).$$
Let $\vec{\mathfrak{a}}: (m) \longrightarrow \bigcup \mathrm{rng}(D)$ be a sequence such that for all $0 \leq l < m$, $\vec{\mathfrak{a}}(l) \in D(j_l)$. Let
$$B= \{\mathfrak{b} \in M_i \mid \exists \vec{\mathfrak{c}}((\vec{\mathfrak{c}}: [m] \longrightarrow V) \land (\vec{\mathfrak{c}}(0)= \mathfrak{b}) \land (\forall l \in (m))(\vec{\mathfrak{c}}(l+1)= \vec{\mathfrak{a}}(l)) \land (\mathcal{D} \models \phi[\vec{\mathfrak{c}}]))\}.$$
Stratified Comprehension ensures that $B$ is a set. Clearly $B \subseteq M_{i+1}$. By Lemma \ref{Th:KeyCompLemma} there exists $\mathfrak{d} \in M_{i+1}$ such that for all $\mathfrak{b} \in M_i$,
$$\begin{array}{ccccc}
\mathcal{D} \models (\mathfrak{b} \in_i \mathfrak \mathfrak{d}) & \textrm{if and only if} & \langle \mathfrak{b}, \mathfrak{d} \rangle \in \mathcal{E} & \textrm{if and only if} & \textrm{there exists } \vec{\mathfrak{c}}: [m] \longrightarrow V \textrm{ with}\\
& & & & \vec{\mathfrak{c}}(0)= \mathfrak{b} \textrm{ and } (\forall l \in (m))(\vec{\mathfrak{c}}(l+1)= \vec{\mathfrak{a}}(l))\\ 
& & & & \textrm{and } \mathcal{D} \models \phi[\vec{\mathfrak{c}}] 
\end{array}$$ 
This shows that Comprehension holds in $\mathcal{D}$.\\
Finally, we need to show that $\mathcal{D} \models (\textrm{Infinity})$. Let
$$R=\{\langle i, j \rangle \mid i, j \in \mathbb{N} \land i < j\} \cup \{\langle i, V \rangle \mid i \in \mathbb{N}\}.$$
The relation $R$ is a BFEXT with $R= \mathrm{seg}_R(V)$, therefore $R \in \Omega$. Let $\mathfrak{a}= [R]$. Note that $\mathfrak{a} \in M_1$. Let $\phi(x, y)$ be the $\mathcal{L}$-formula
$$(x= \langle z_1, z_2 \rangle) \land (z_1, z_2 \in y) \land (z_2= z_1 \cup \{z_1\}).$$
Let
$$B= \{ \mathfrak{b} \in M_0 \mid \exists \vec{\mathfrak{c}} ((\vec{\mathfrak{c}}: (2) \longrightarrow V) \land (\vec{\mathfrak{c}}(0)= \mathfrak{b}) \land (\vec{\mathfrak{c}}(1)= \mathfrak{a}) \land (\langle M_1, \mathcal{E} \rangle \models \phi[\vec{\mathfrak{c}}]))\}.$$
Stratified Comprehension ensures that $B$ is a set. We also have that $B \subseteq M_3$. Using Lemma \ref{Th:KeyCompLemma} we can find $\mathfrak{d} \in M_3$ such that for all $\mathfrak{b} \in M_2$,
$$\langle \mathfrak{b}, \mathfrak{d} \rangle \in \mathcal{E} \textrm{ if and only if } \mathfrak{b} \in B.$$
The point $\mathfrak{d}$ in $\mathcal{D}$ is an injective function that witnesses that $\mathfrak{a}$ is Dedekind infinite.\\
This shows that $\mathcal{M} \models \mathrm{Con}(\mathrm{TSTI}_{n+1})$. Since $k < n+1$, $\mathcal{M} \models \mathrm{Con}(\mathrm{TSTI}_k)$, which is a contradiction.               
\Square
\end{proof}

Since $\mathcal{M}$ was an arbitrary model of $\mathrm{NFU}^{-\mathrm{AC}}+\mathrm{AxCount}_\geq+ \neg\mathrm{AxCount}$ this proves:

\begin{Theorems1} \label{Th:MainTheorem}
$\mathrm{NFU}^{-\mathrm{AC}}+\mathrm{AxCount}_\geq \vdash \mathrm{Con}(\mathrm{TSTI})$. \Square
\end{Theorems1}

\noindent Combined with Theorem \ref{Th:ConsistencyStregthOfNFU} this shows that the theory $\mathrm{NFU}^{-\mathrm{AC}}+\mathrm{AxCount}_\geq$ has strictly stronger consistency strength than the theory $\mathrm{NFU}^{-\mathrm{AC}}$.

\begin{Coroll1} \label{Th:AxCountGEQStrictlyStrongerOverNFU}
$\mathrm{NFU}^{-\mathrm{AC}}+\mathrm{AxCount}_\geq \vdash \mathrm{Con}(\mathrm{NFU}^{-\mathrm{AC}})$. \Square
\end{Coroll1}

Again, Theorem \ref{Th:StrongConsistencyOfChoiceWithNFU} allows any occurrence of ``$\mathrm{NFU}^{-\mathrm{AC}}$" in Corollary \ref{Th:AxCountGEQStrictlyStrongerOverNFU} to be replaced with ``$\mathrm{NFU}$". The following from \cite{mck15} still remains open:

\begin{Quest1}
What is the exact consistency strength of $\mathrm{NFU}+\mathrm{AxCount}_\geq$ relative to a subsystem of $\mathrm{ZFC}$?
\end{Quest1}

Since the theory $\mathrm{NF}$ can be viewed as an extension of the theory $\mathrm{NFU}^{-\mathrm{AC}}$, Theorem \ref{Th:MainTheorem} also yields:

\begin{Coroll1}
$\mathrm{NF}+\mathrm{AxCount}_\geq \vdash \mathrm{Con}(\mathrm{TSTI})$. \Square
\end{Coroll1}

\bibliographystyle{alpha}
\bibliography{.}          

\begin{thebibliography}{9}

\bibitem[Ena]{ena06} Enayat, Ali. ``From bounded arithmetic to second order arithmetic via automorphisms". \emph{Logic in Tehran}. Lecture Notes in Logic. Vol. 26. Association for Symbolic Logic. 2006.

\bibitem[For77]{for77} Forster, Thomas E. \emph{N.F.} Ph.D. Thesis. University of Cambridge, United Kingdom. 1977.

\bibitem[For95]{for95} Forster, Thomas E. \emph{Set Theory with a Universal Set: Exploring an Untyped Universe}. Oxford Logic Guides no 31. 1995.

\bibitem[FH]{fh09} Forster, Thomas E. and Holmes, M. Randall. ``Permutation methods in $\mathrm{NF}$ and $\mathrm{NFU}$". \emph{Proceedings of the $70^\textrm{th}$ anniversary $\mathrm{NF}$ meeting in Cambridge}. Edited by M. Crabbe and T. E. Forster. Cahiers du Centre de Logique. Vol. 16. Academia, Louvain-la-Neuve. 2009. pp 33-76

\bibitem[Hin]{hin75} Hinnion, Roland. \emph{Sur la th\'{e}orie des ensembles de Quine}. PhD Thesis. ULB, Brussels. 1975. Translated by Thomas Forster. Available online from http://www.logic-center.be/Publications/Bibliotheque/hinnionthesis.pdf. 2009.

\bibitem[Hol98]{hol98} Holmes, M. Randall. \emph{Elementary Set Theory with a Universal Set}. Cahiers du Centre de logique. Vol. 10. Academia, Louvain-la-Neuve. 1998.

\bibitem[Hol01]{hol01} Holmes, M. Randall. ``Strong Axioms of Infinity in $\mathrm{NFU}$". \emph{The Journal of Symbolic Logic}. Vol. 66. 2001. pp 87-116

\bibitem[Jen]{jen69} Jensen, Ronald B. ``On the Consistency of a Slight (?) Modification of Quine's New Foundations". \emph{Synthese}. Vol. 19. 1969. pp 250-263

\bibitem[M]{mck15} McKenzie, Zachiri. ``Automorphisms of models of set theory and extensions of $\mathrm{NFU}$". \emph{Annals of Pure and Applied Logic}. Vol. 166. 2015. pp 601-638

\bibitem[Mat]{mat01} Mathias, Adrian R. D. ``The strength of Mac Lane set theory". \emph{Annals of Pure and Applied Logic}. Vol. 110. 2001. pp 107-234

\bibitem[Ore]{ore64} Orey, Steven. ``New Foundations and the axiom of counting". \emph{Duke Mathematical Journal}. Vol. 31. 1964. pp 655-660

\bibitem[Qui37]{qui37} Quine, Willard v. O. ``New foundations for mathematical logic". \emph{American Mathematical Monthly}. Vol. 44. 1937. pp 70-80

\bibitem[Qui45]{qui45} Quine, Willard v. O. ``On ordered pairs". \emph{Journal of Symbolic Logic}. Vol. 10. 1945. pp 95-96

\bibitem[Ros]{ros78} Rosser, J. Barkley. \emph{Logic for mathematicians}. McGraw-Hill, reprinted (with appendices) by Chelsea, New York. 1978.

\bibitem[RW]{rw08} Russell, Bertrand A. W. and Whitehead, Alfred N. \emph{Principia Mathematica}. Cambridge University Press. 1908.

\bibitem[Spe]{spe53} Specker, Ernst P. ``The Axiom of Choice in Quine's ``New Foundations for Mathematical Logic"". \emph{Proceedings of the National Academy of Sciences, U.S.A.} Vol. 29. 1953. pp 366-368

\bibitem[Sol]{solXX} Solovay, Robert. ``The consistency strength of $\mathrm{NFUB}$". Preprint. Available online from http://arxiv.org/.

\end{thebibliography}

\end{document}